\documentclass[11pt,twoside]{article}
\usepackage{latexsym,amsmath}
\usepackage{indentfirst}
\usepackage{graphicx}

\numberwithin{equation}{section}
\newtheorem{Prop}{\bf Proposition}[section]

\newtheorem{defn}{\bf Definition}[section]
\newtheorem{Rem}{\bf Remark}[section]
\newtheorem{Ex}{\bf Example}[section]

\begin{document}
\def \b{\Box}

\begin{center}
{\Large {\bf The $q-$fractional  Euler top system with one control\\[0.1cm]
on a fractional Leibniz algebroid }}

\end{center}

\begin{center}
{\bf Gheorghe IVAN}
\end{center}

\setcounter{page}{1}
\pagestyle{myheadings}

{\small {\bf Abstract}. We determine the fractional almost Poisson
realizations for fractional Euler top system with one control. These
realizations allow us to introduce a fractional Leibniz algebroid
structure  on ${\bf R}^{3}$ and also to define the $q-$fractional
Euler top system with one control on a fractional Leibniz
algebroid.
  Finally, the numerical integration of them are discussed.
  {}

\section {Introduction}
\smallskip

Fractional calculus is an interdisciplinary area of applied mathematics with many applications in several fields of science and engineering. For example, the fractional calculus has been applied to variational problems in: fluid mechanics,  geometrical mechanics, chaotic dynamics,  quantum  mechanics, nanotechnology \cite{podl, agra, klag, nabu, kilb, elna,  bzte, igom, gi22}. The fractional derivatives and integrals describe with great precision the dynamic behavior of various systems of differential equations \cite {tara, imgo, migo, yama, opmi, enls, ghiv, mist, mlag}.

The Lie algebroids have been used to  geometric formulation of Hamiltonian mechanics.
In Hamiltonian mechanics the approach is to consider  the cotangent bundle. The cotangent bundle  forms a phase space consisting of the configuration and momentum space.  A state is a point in the phase space and  an observable is given by Hamiltonian function. This function is used to generate the sections of the bundle. The evolution is given by a bundle morphism (named,  anchor) and is defined by Hamilton’s equations.
The concept of Leibniz algebroid is a wakened version of a Lie algebroid, where the bilinear operation on sections of the vector bundle is not necessarily skew-symmetric and the anchor is a homomorphism of Leibniz algebra.

The Lie algebroids and Leibniz algebroids have proven to be powerful tools for geometric formulation of the Hamiltonian mechanics \cite{mana, bagr, demi,nabu1, ezuc, mart, pope}. Also, they have been used  in the investigation of many fractional dynamical systems \cite{gmio, imod, pagi, chdo, igmp}.

In this paper we refer to a special class of dynamical
systems. This class is formed by a family of differential equations on ${\bf R}^{3}$ which depends by a triple
of real parameters $(\alpha_{1}, \alpha_{2}, \alpha_{3}) $ and one
control parameter $k\in {\bf R},$  called the Euler top system with one control.
It is described by the following nonlinear differential equations on $~{\bf R}^{3}:$\\[-0.4cm]
\begin{equation}
 \dot{x}^{1}  = \alpha_{1} x^{2} x^{3} ,~~~\dot{x}^{2}  = \alpha_{2}  x^{1} x^{3} - k x^{1},~~~\dot{x^{3}}  = \alpha_{3} x^{1} x^{2} ,\label{(1.1)}
\end{equation}
where $~x^{1}, x^{2}, x^{3}~$ are state variables and $ \alpha_{i}, k\in{\bf R} $ for $ i=\overline{1,3}$ such that $
 \alpha_{1} \alpha_{2} \alpha_{3}\neq 0 .$

The family of dynamical systems $(1.1)$ generalizes for example: the Lagrange system \cite{takh}, a particular case of  the Rabinovich system  which  is determined by the nonlinear terms \cite{chdo} and  the general Euler top system \cite{ivan, ivmi}.

The goal of our work is to study the fractional-order model associated to Euler top system with one control $(1.1).$

The paper is structured as follows. In Section 2, fractional Euler
top system with one control  $(2.4)$ is presented. We prove that the
fractional top system $(2.4)$ has two fractional almost Poisson realizations
(Proposition 2.1). In Section 3, the Proposition 2.1 is used to
give the construction of $q-$fractional Euler top system with control on a fractional Leibniz algebroid. In Section 4,
the numerical integration $q-$fractional Euler top system with one control $(3.8)$ is discussed. As example, the numerical simulation for the $q-$fractional Lagrange system with one
control $(3.9)$  is given.

\section{ Fractional almost Poisson realizations of fractional Euler top system with one control}

We start by recalling of the notion fractional Leibniz structure
\cite{ivop}.
\markboth{Gheorghe IVAN}{The $q-$fractional Euler top system with one control on a ...}

Let $ f\in C^{\infty}( \textbf{R}) $ and $ q \in \textbf{R}, q
> 0. $ The $ q-$order Caputo differential operator
\cite{bagr, die1} is described by $~D_{t}^{q}f(t) = I^{m -
}f^{(m)}(t),~q > 0,$ where $~f^{(m)}(t)$  is the $ m-$order
derivative of the function $ f,~m \in \textbf {N}^{\ast}$ is an
integer such that $ m-1 \leq q \leq m $ and $ I^{q} $ is the $
q-$order Riemann-Liouville integral
operator expressed by:\\[-0.1cm]
\[
I^{q}f(t) =\displaystyle\frac{1}{\Gamma(q)}\int_{0}^{t}{(t-s)^{q
-1}}f(s)ds,~~~q > 0,
\]
where $ \Gamma $ is the Euler Gamma function. If $ q =1$, then $ D_{t}^{q}f(t) = df / dt.~$ We suppose $~q\in (0,1).$

 Let $ M $ be a $ n-$ dimensional smooth manifold, $ U \subset M $ a local chart, $ ( x^{i}),
 i=\overline{1,n} $ a  system of coordinates on $ U $ and $ D_{x^{i}}^{q}f $ the Caputo partial derivatives  for $ f\in
C^{\infty}(U). $

We denote by  $ \mathcal{X}^{q}(U)$ be the module
of the fractional vector fields generated by the operators $
D_{x^{i}}^{q}, i = \overline{1,n}. $ A fractional vector field $
{\overset{q}{X}}\in \mathcal{X}^{\alpha}(U) $ \cite{ivop}, has the
form:\\[-0.2cm]
\[
{\overset{q}{X}} = {\overset{q}{X}}^{i} D_{x^{i}}^{q},~
{\overset{q}{X}}^{i}\in C^{\infty}(U),~i=\overline{1,n}.
\]
Let $ \mathcal{D}^{q}(U) $ the module
 generated by $1-$forms $~ d( x^{i})^{q}, i=\overline{1,n} $ on $ U.~$
The  fractional exterior derivative $ d^{q}: C^{\infty}(U)
\rightarrow \mathcal{D}^{q}(U)$  is expressed by:\\[-0.2cm]
\[
d^{q}(f)= d(x^{i})^{q} D_{x^{i}}^{q}(f),~f\in C^{\infty}(U).
\]
Let $ {\overset{q}{P}}\in \mathcal{X}^{q}(M)\times
\mathcal{X}^{q}(M) $  be a  fractional $ 2- $ tensor field and $
d^{q} f, d^{q}g \in \mathcal{D}^{q}(M) $. The bilinear map $
[\cdot,\cdot]^{q} : C^{\infty}(M) \times C^{\infty}(M) \rightarrow
C^{\infty}(M)$  defined by:\\[-0.2cm]
\[
[f,h]^{q} = {\overset{q}{P}}( d^{q}f,d^{q}h ), ~~~ (\forall )
f,h\in C^{\infty}(M),
\]
 is called the \textit{ fractional Leibniz bracket}.

The \textit{fractional Leibniz vector field  $
{\overset{q}{X}}_{h} $ associated to $ h\in C^{\infty}(M),$} is defined by:\\[-0.2cm]
\[
{\overset{q}{X}}_{h}(f)= [f,h]^{q}, ~ (\forall) f\in
C^{\infty}(M),
\]
 Locally, if $~{\overset{q}{P}} $ is a skew-symmetric
fractional $ 2-$tensor field on $ M $ generated by the matrix $
{\overset{q}{P}}= ( {\overset{q}{P}}^{ij}),$ then the dynamical
system associated to $ {\overset{q}{X}}_{h} $ is given
by:\\[-0.2cm]
\[
D_{t}^{q} x^{i}(t)= [x^{i}(t),h(t)]^{q}, ~~~\textrm{ where}~~~
[x^{i},h]^{q} = {\overset{q}{P}}^{ij}\cdot
D_{x^{j}}^{q}h,\\[-0.2cm]
\]
or equivalently\\[-0.4cm]
\begin{equation}
D_{t}^{q} x^{i} = {\overset{q}{P}}^{ij}\cdot
D_{x^{j}}^{q}h.\label{(2.1)}\\[-0.2cm]
\end{equation}

The system $(2.1)$ is called the {\it fractional Leibniz system}
associated to matrix $ {\overset{\alpha}{P}} $ with the
Hamiltonian $ h\in C^{\infty}(M) $.

Let now $ [\cdot,\cdot ]^{q} $ be a fractional Leibniz structure
on $ {\bf R}^{n}$ with the associated matrix $ {\overset{q}{P}} =
( {\overset{q}{P}}^{ij}), $  i.e. $ {\overset{q}{P}}^{ij} = [
x^{i}, x^{j}]^{q},~i,j=\overline{1,n}.$
\begin{defn}
{\rm (\cite{ivop}) Let be a fractional dynamical system on ${\bf
R}^{n}$ given by:\\[-0.2cm]
\begin{equation}
D_{t}^{q}x^{i}(t) = f^{i}( x^{1}(t),\ldots,x^{n}(t)),~~~f^{i}\in
C^{\infty}({\bf R}^{n}, {\bf R}),~i=\overline{1,n}.\label{(2.2)}
\end{equation}

We say that the fractional dynamics $ (2.2) $ has a {\it
fractional Leibniz realization}, if there exist a Leibniz
structure $ [\cdot,\cdot ]^{q} $ on $ {\bf R}^{n}$ generated by
the matrix $ {\overset{q}{P}} = ( {\overset{q}{P}}^{ij}) $ and a
Hamiltonian $ {\overset{q}{H}}\in C^{\infty}({\bf R}^{n}, {\bf
R})$ such that $(2.2) $ can be written in the following
form:\\[-0.2cm]
\begin{equation}
D_{t}^{q}x(t) = {\overset{q}{P}}(x(t))\cdot \nabla^{q}
{\overset{q}{H}}(x(t)),\label{(2.3)}\\[-0.1cm]
\end{equation}
where $ x(t) = ( x^{1}(t), \ldots,x^{n}(t))^{T}~$ and $~\nabla^{q}
{\overset{q}{H}}= (
D_{x^{1}}^{q}{\overset{q}{H}},\ldots,D_{x^{n}}^{q}{\overset{q}{H}})^{T}$.}\hfill$\Box$
\end{defn}

 A fractional Leibniz realization is denoted by $( {\bf R}^{n},
[\cdot, \cdot]^{q}, {\overset{q}{H}} ) $ or $ ({\bf R}^{n},
{\overset{q}{P}}, {\overset{q}{H}}).$

If the matrix $  {\overset{q}{P}} = ( {\overset{q}{P}}^{ij}) $ is
skew-symmetric, then $ ({\bf R}^{n}, {\overset{q}{P}},
{\overset{q}{H}} ) $ is called the {\it fractional almost Poisson
realization} of the system $(2.2).$

A {\it fractional Casimir} of the configuration $( {\bf R}^{n},
[\cdot, \cdot]^{q}, {\overset{q}{H}} ) $ is a function
${\overset{q}{C}}\in C^{\infty}({\bf R}^{n}, {\bf R})$ such that
$[ {\overset{q}{C}}, f]^{q} =0 $ for all $ f\in C^{\infty}({\bf
R}^{n}, {\bf R}),$ that is $~{\overset{q}{P}}(x(t))\cdot
\nabla^{q} {\overset{q}{C}}(x(t)) = 0.$

The {\it fractional Euler top  system with parameters} is
defined by the following set of fractional differential equations:\\[-0.2cm]
\begin{equation}
\left\{ \begin{array} {lcl}
 D_{t}^{q}{x}^{1}(t) & = & \alpha_{1}
 x^{2}(t) x^{3}(t)  \\[0.1cm]
 D_{t}^{q}{x}^{2}(t) & = & \alpha_{2}
 x^{1}(t) x^{3}(t) - k x^{1}(t) ,~~~ q \in (0,1),\\[0.1cm]
  D_{t}^{q}{x}^{3}(t) & = & \alpha_{3} x^{1}(t)x^{2}(t) \label{(2.4)}
  \end{array}\right.
\end{equation}
 where $ \alpha_{i}, k\in{\bf R} $ for $ i=\overline{1,3}$ such that $
 \alpha_{1} \alpha_{2} \alpha_{3}\neq 0 .$

For $k\neq 0,$ the fractional dynamics  $ (2.4)$ is called the
{\it fractional Euler top system with one control}. If $
k=0, $ then $ (2.4)$ is the {\it fractional Euler top system}, \cite{migi}.

For $ q=1,~(2.4) $ becomes the {\it Euler top system with
parameters} given by\\[-0.2cm]
\begin{equation}
 \dot{x}^{1}  =  \alpha_{1}
 x^{2} x^{3},~~
 \dot{x}^{2}  = \alpha_{2}
 x^{1} x^{3} - k x^{1},~~
  \dot{x}^{3}=  \alpha_{3} x^{1}x^{2}, \label{(2.5)}
\end{equation}
 where $ \alpha_{i}, k \in{\bf R} $ for $ i=\overline{1,3}$ such that $
 \alpha_{1} \alpha_{2} \alpha_{3}\neq 0. $

 For $k\neq 0,$ the dynamics $ (2.5)$ is called the
{\it Euler top system with one control}. If $ k=0, $ then $
(2.5)$ is the {\it Euler top system} \cite{igim, ivan}.

\begin{Ex}
{\rm $(i)~$ For $ \alpha = ( 1 , 1, 1) $
and $ k\neq 0,$ the system  $(2.4)$ becomes:\\[-0.2cm]
\begin{equation}
 D_{t}^{q} x^{1}  =  x^{2}x^{3},~~~
 D_{t}^{q} x^{2}  =  x^{1}x^{3} -  k x^{1},~~~
 D_{t}^{q} x^{3}  =  x^{1}x^{2}.\label{(2.6)}
\end{equation}

The system $(2.6) $ is called the {\it fractional Lagrange system
with one control}. If $ k = 0 $ and $q=1,$ the system
$(2.6)$ is the {\it Lagrange system} in the theory of static
$SU(2)-$monopoles \cite{takh}.

$(ii)~$ For $ \alpha = ( 1 , - 1, 1)$
and $ k\neq 0,$ the system  $(2.4)$ becomes:\\[-0.2cm]
\begin{equation}
 D_{t}^{q} x^{1}  =  x^{2} x^{3},~~~
 D_{t}^{q} x^{2}  =  - x^{1}x^{3} - k x^{1},~~~
 D_{t}^{q}x^{3}  =  x^{1}x^{2}.\label{(2.7)}
\end{equation}

The system $(2.7) $ is called the {\it fractional Rabinovich
system with one control}. For $ k = 0 $ and $q=1,$ the
system $(2.7)$ is the {\it Rabinovich system} \cite{chdo}.}
\hfill$\Box$
\end{Ex}
\begin{Prop}
Let be $ q\in (0,1]. $ The fractional Euler top  system with
parameters $(2.4)$ has the following fractional almost Poisson
realizations:

$(i)~~( {\bf R}^{3}, P_{1}^{\alpha},
{\overset{q}{H}}_{1}^{\alpha}) $ with the fractional Casimir $
{\overset{q}{C}}_{1}^{\alpha},$ where
\begin{equation} P_{1}^{\alpha}(x,k) = \left ( \begin{array}{ccc}
0 & k & -\alpha_{3} x^{2}\\[0.2cm]
-k & 0 & -\displaystyle\frac{\alpha_{2}\alpha_{3}}{\alpha_{1}} x^{1}\\[0.3cm]
\alpha_{3} x^{2} & \displaystyle\frac{\alpha_{2}\alpha_{3}}{\alpha_{1}} x^{1} & 0\\
\end{array}\right ),\label{(2.8)}
\end{equation}
\begin{equation}
{\overset{q}{H}}_{1}^{\alpha}(x) =
\displaystyle\frac{1}{\Gamma(q+2)}[(x^{1})^{1+q}
-\displaystyle\frac{\alpha_{1}}{\alpha_{3}} (x^{3})^{1+q}],
\label{(2.9)}
\end{equation}
\begin{equation}
{\overset{q}{C}}_{1}^{\alpha}(x,k) =
\displaystyle\frac{1}{\Gamma(q+2)}[\alpha_{2}(x^{1})^{1+q} -
\alpha_{1}(x^{2})^{1+q}] - \displaystyle\frac{k
\alpha_{1}}{\Gamma(q+1)\alpha_{3}} (x^{3})^{q}. \label{(2.10)}
\end{equation}

$(ii)~~( {\bf R}^{3}, P_{2}^{\alpha},
{\overset{q}{H}}_{2}^{\alpha}) $ with the fractional Casimir $
{\overset{q}{C}}_{2}^{\alpha},$ where
\begin{equation} P_{2}^{\alpha}(x) = \left ( \begin{array}{ccc}
0 &
x^{3} &  0\\[0.2cm]
- x^{3} & 0 & -\displaystyle\frac{\alpha_{3}}{\alpha_{1}} x^{1}\\[0.2cm]
0 &  \displaystyle\frac{\alpha_{3}}{\alpha_{1}} x^{1} & 0\\[0.2cm]
\end{array}\right ),\label{(2.11)}
\end{equation}
\begin{equation}
{\overset{q}{H}}_{2}^{\alpha}(x,k) =
\displaystyle\frac{1}{\Gamma(q+2)}[-\alpha_{2} (x^{1})^{1+q} +
\alpha_{1}(x^{2})^{1+q}]+
\displaystyle\frac{k\alpha_{1}}{\Gamma(q+1)\alpha_{3}}(x^{3})^{q},
\label{(2.12)}
\end{equation}
\begin{equation}
{\overset{q}{C}}_{2}^{\alpha}(x) =
\displaystyle\frac{1}{\Gamma(q+2)}[\alpha_{3}(x^{1})^{1+q} -
\alpha_{1} (x^{3})^{1+q}]. \label{(2.13)}
\end{equation}
\end{Prop}
{\it Proof.} $(i)~$ Using the relations $~D_{x^{i}}^{q}
(x^{j})^{m} = \delta_{i}^{j}\frac{\Gamma(1+ m)}{\Gamma (1+m-q)}
(x^{j})^{m - q}, i,j=\overline{1,3} $ follows:\\[0.1cm]
 $~D_{x^{1}}^{q}{\overset{q}{H}}_{1}^{\alpha}= x^{1},~
D_{x^{2}}^{q}{\overset{q}{H}}_{1}^{\alpha} = 0
,~D_{x^{3}}^{q}{\overset{q}{H}}_{1}^{\alpha}=
-\displaystyle\frac{\alpha_{1}}{\alpha_{3}}x^{3},~
D_{x^{1}}^{q}{\overset{q}{C}}_{1}^{\alpha}= \alpha_{2}
x^{1},~D_{x^{2}}^{q}{\overset{q}{C}}_{1}^{\alpha}= - \alpha_{1}
x^{2},$\\[0.1cm]
$D_{x^{3}}^{q}{\overset{q}{C}}_{1}^{\alpha}= -\displaystyle\frac{k
\alpha_{1}}{\alpha_{3}}.$ Then: $~~P_{1}^{\alpha}\cdot \nabla
{\overset{q}{H}}_{1}^{\alpha}  =$\\
\[
= \left ( \begin{array}{ccc} 0 &
k & -\alpha_{3} x^{2} \\[0.2cm]
-k & 0 & -\displaystyle\frac{\alpha_{2}\alpha_{3}}{\alpha_{1}} x^{1}\\[0.2cm]
\alpha_{3} x^{2} &  \displaystyle\frac{\alpha_{2}\alpha_{3}}{\alpha_{1}} x^{1} & 0\\
\end{array}\right ) \left ( \begin{array}{c}
x^{1} \\[0.2cm]
 0\\[0.2cm]
-\displaystyle\frac{\alpha_{1}}{\alpha_{3}} x^{3} \\
\end{array}\right ) = \left (\begin{array}{c}
\alpha_{1} x^{2} x^{3}\\[0.2cm]
\alpha_{2} x^{1}x^{3} - k x^{1} \\[0.2cm]
\alpha_{3} x^{1}x^{2}\\
\end{array}\right )=\left (\begin{array}{c}
D_{t}^{q} x^{1} \\[0.2cm]
D_{t}^{q}x^{2}\\[0.2cm]
D_{t}^{q} x^{3}\\
\end{array}\right ).
\]

Hence  $ D_{t}^{q} x(t)= P_{1}^{\alpha}(x(t)) \cdot \nabla
{\overset{q}{H}}_{1}^{\alpha}(x(t)) $ and $ (2.4) $ is a
fractional almost Poisson system. Also, $ {\overset{q}{C}}_{1}^{\alpha} $ is a fractional Casimir,
since:\\[-0.2cm]
 \[
 P_{1}^{\alpha} \cdot \nabla {\overset{q}{C}}_{1}^{\alpha} =
\left ( \begin{array}{ccc} 0 &
k & -\alpha_{3} x^{2} \\[0.2cm]
-k & 0 & -\displaystyle\frac{\alpha_{2}\alpha_{3}}{\alpha_{1}} x^{1}\\[0.2cm]
\alpha_{3} x^{2} &  \displaystyle\frac{\alpha_{2}\alpha_{3}}{\alpha_{1}} x^{1} & 0\\
\end{array}\right ) \left ( \begin{array}{c}
\alpha_{2} x^{1} \\[0.2cm]
 -\alpha_{1} x^{2}\\[0.2cm]
-\displaystyle\frac{k\alpha_{1}}{\alpha_{3}}\\
\end{array}\right )= \left (\begin{array}{c}
0 \\[0.2cm]
0\\[0.2cm]
0\\
\end{array}\right ).
\]

$(ii)~$ We have $ D_{x^{1}}^{q}{\overset{q}{H}}_{2}^{\alpha}=
-\alpha_{2} x^{1},~ D_{x^{2}}^{q}{\overset{q}{H}}_{2}^{\alpha} =
\alpha_{1} x^{2}
,~D_{x^{3}}^{q}{\overset{q}{H}}_{2}^{\alpha}=\displaystyle\frac{k
\alpha_{1}}{\alpha_{3}},~
D_{x^{1}}^{q}{\overset{q}{C}}_{2}^{\alpha}= \alpha_{3}
x^{1},~D_{x^{2}}^{q}{\overset{q}{C}}_{2}^{\alpha}=
0,~~D_{x^{3}}^{q}{\overset{q}{C}}_{2}^{\alpha}= - \alpha_{1}
x^{3}.$ Similarly, $~D_{t}^{q} x = P_{2}^{\alpha} \cdot \nabla
{\overset{q}{H}}_{2}^{\alpha}~$ and $~ P_{2}^{\alpha} \cdot \nabla
{\overset{q}{C}}_{2}^{\alpha} = 0.~\hfill\Box$

\begin{Rem}
{\rm If in Proposition 2.1 we replace $ \alpha= (\alpha_{1},
\alpha_{2}, \alpha_{3}) $  with $ \alpha = (1,1,1)$ (resp. $
\alpha = (1,-1,1)$) we obtain two fractional almost Poisson
realizations for the fractional Lagrange (resp. Rabinovich) system
with parameter $(2.6)$ (resp. $(2.7)$)}. \hfill$\Box$
\end{Rem}

\section{ Fractional Euler top system with one control on a fractional Leibniz algebroid}

A {\t fractional Leibniz algebroid structure} \cite{imod} on
a vector bundle $ \pi : E \rightarrow M $ is
 given by a bracket $ [ \cdot, \cdot ]^{q} $ on the space of sections $ Sec(\pi) $ and two vector bundle morphisms
$ {\overset{q}{\rho}}_{1}, {\overset{q}{\rho}}_{2}: E \rightarrow
T^{q}M $ (called the {\it left} and {\it right fractional anchor})
such that for all $ \sigma_{1}, \sigma_{2}\in Sec(\pi)$ and $ f, g
\in C^{\infty}(M)$ we have:\\[-0.2cm]
\[
\left \{ \begin{array}{l} [e_{u}, e_{v}]^{q}
=C_{uv}^{w}e_{w},\\[0.2cm]
 [ f \sigma_{1}, g \sigma_{2} ]^{q}
=f {\overset{q}{\rho}}_{1}(\sigma_{1})(g)\sigma_{2} - g
{\overset{q}{\rho}}_{2}(\sigma_{2})(f) \sigma_{1} + f g [
\sigma_{1}, \sigma_{2}]^{q}.
\end{array}\right.
\]

 A vector bundle  $ \pi : E\rightarrow M $ endowed with a fractional Leibniz algebroid structure
 $ ( [\cdot, \cdot ]^{q}, {\overset{q}{\rho}}_{1}, {\overset{q}{\rho}}_{2} ) $ on $ E $ , is called {\it
  fractional Leibniz algebroid} over $ M $ and it is denoted by $(E, [\cdot, \cdot ]^{q}, {\overset{q}{\rho}}_{1}, {\overset{q}{\rho}}_{2}).$

In the paper \cite{imod} has shown that there exists a bijective
correspondence between the fractional Leibniz algebroid structures
on a vector bundle $ E $ and the linear fractional $ 2-$tensor
fields on the dual vector bundle $ E^{\ast}.$

More precisely, a linear fractional $ 2-$tensor field $
{\overset{q}{\Lambda}} $ on the dual vector bundle $ \pi^{\ast} :
E^{\ast} \to M $ defines a fractional Leibniz structure on the
vector bundle $\pi : E \to M.$

The bracket $ [\cdot, \cdot ]_{{\overset{q}{\Lambda}}} $ is
defined by:\\[-0.3cm]
 \[
 [ F, G ]_{{\overset{q}{\Lambda}}^{r}} = {\overset{q}{\Lambda}}^{r}( d^{q r}F, d^{q r} G ), ~(\forall) ~ F,G\in
 C^{\infty}(E^{\ast}),~~~~~\hbox{where}\\[-0.1cm]
 \]
 \[
 d^{q r}F = d(x^{i})^{q} D_{x^{i}}^{q}H + d(\xi_{u})^{r}D_{\xi_{u}}^{r}F= d^{q}(F)+d^{r}(F).
 \]

If $ ( x^{i} ),~ ( x^{i}, y^{u} ) $ resp., $ ( x^{i}, \xi_{u} ) $
for $ i=\overline{1,n},~ u =\overline{1,m} $ are coordinates on $
M, E $ resp. $ E^{\ast} $, then the linear fractional tensor
${\overset{q}{\Lambda}}^{r} $ on $ E^{\ast} $ has the form:\\[-0.3cm]
\[
{\overset{q}{\Lambda}}^{r} = C_{uv}^{w}\xi_{w}
D_{\xi_{u}}^{r}\otimes D_{\xi_{v}}^{r} + {\overset{q}{\rho}}_{1
u}^{i} D_{\xi_{u}}^{r}\otimes D_{x^{i}}^{q}  -
{\overset{q}{\rho}}_{2u}^{i}D_{x^{i}}^{q}\otimes D_{\xi_{u}}^{r}.
\]

The fractional system associated to vector field $
{\overset{q}{X}}_{H}^{r} $ with $ H\in C^{\infty}(E^{\ast}) $
given by:\\[-0.2cm]
\begin{equation}
{\overset{q}{X}}_{H}^{r}(F) = {\overset{q}{\Lambda}}^{r}(d^{q r}F,
d^{q r}H),~\hbox { for all }~ F\in C^{\infty}(E^{\ast})\label
{(3.1)}
\end{equation}
is called a {\it fractional dynamical system on the fractional Leibniz algebroid} $(E, [\cdot,
\cdot ]^{q}, {\overset{q}{\rho}}_{1}, {\overset{q}{\rho}}_{2} ),$ \cite{ivop}.

 Locally, the fractional dynamical system $ (3.1) $
 reads:\\[-0.2cm]
\begin{equation}
\left\{\begin{array}{l}
 D_{t}^{q} \xi_{u} = [ \xi_{u},
H]_{{\overset{q}{\Lambda}}^{r}}=C_{uv}^{w}\xi_{w}D_{\xi_{v}}^{r}
H + {\overset{q}{\rho}}_{1u}^{i}D_{x^{i}}^{q} H\\
D_{t}^{q}x^{i} = [ x^{i}, H]_{{\overset{q}{\Lambda}}^{r}}=
-{\overset{q}{\rho}}_{2u}^{i}D_{\xi_{u}}^{r} H \\
\end{array}\right.
. \label {(3.2)}\\[-0.1cm]
\end{equation}

{\bf In the following we suppose that $ r=q.$}

 If $~\widetilde{P} = (C_{uv}^{w}\xi_{w}),~
\overset{q}{\rho}_{1}=({\overset{q}{\rho}}_{1u}^{i})~$ and $~
\overset{q}{\rho}_{2}=({\overset{q}{\rho}}_{2u}^{i}),$ for
$~i=\overline{1,3},~u,v,w=\overline{1,3},$ then the fractional
system $ (3.2) $ (with condition $ r=q $) can be written in the
matrix form:\\[-0.3cm]
\begin{equation}
\left (\begin{array}{c}
D_{t}^{q}{\xi}_{1}\\[0.2cm]
D_{t}^{q}{\xi}_{2}\\[0.2cm]
D_{t}^{q}{\xi}_{3}\\
\end{array}\right )= \widetilde{P} \left (\begin{array}{c}
D_{\xi_{1}}^{q}H\\[0.2cm]
D_{\xi_{2}}^{q}H\\[0.2cm]
D_{\xi_{3}}^{q}H\\
\end{array}\right ) + \overset{q}{\rho}_{1} \left (\begin{array}{c}
D_{x^{1}}^{q}H\\[0.2cm]
D_{x^{2}}^{q}H\\[0.2cm]
D_{x^{3}}^{q}H\\
\end{array}\right ),~ \left( \begin{array}{c}
D_{t}^{q} x^{1}\\[0.2cm]
D_{t}^{q} x^{2}\\[0.2cm]
D_{t}^{q} x^{3}\\
\end{array}\right )= - \overset{q}{\rho}_{2} \left (\begin{array}{c}
D_{\xi_{1}}^{q}H\\[0.2cm]
D_{\xi_{2}}^{q}H\\[0.2cm]
D_{\xi_{1}}^{q}H\\
\end{array}\right ).\label{(3.3)}
\end{equation}

As example we give the {\it fractional Euler top system with
one control on a fractional Leibniz algebroid}.

 Let the vector bundle $ \pi : E = {\bf R}^{3}\times {\bf
R}^{3}\to {\bf R}^{3} $ and
 $ \pi^{*} : E^{*}={\bf R}^{3}\times ({\bf R}^{3})^{*}\to {\bf R}^{3} $ its dual.
 We consider on $ E^{*} $ the
 linear $ 2-$tensor field $ \overset{q}{\Lambda}, $ defined by the matrix $ \widetilde{P},$ the fractional anchors
 $~\overset{q}{\rho}_{1}, \overset{q}{\rho}_{2}: {\bf R}^{3}\times {\bf R}^{3} \rightarrow
T({\bf R}^{3})~$ and the function $ H $ defined by:
\begin{equation}
\widetilde{P}=P(x,\xi,k) = \left ( \begin{array}{ccc}
0 & k & -\alpha_{3} x^{2}\xi_{2}\\[0.2cm]
-k & 0 & -\displaystyle\frac{\alpha_{2}\alpha_{3}}{\alpha_{1}} x^{1} \xi_{1}\\[0.3cm]
\alpha_{3} x^{2}\xi_{2} & \displaystyle\frac{\alpha_{2}\alpha_{3}}{\alpha_{1}} x^{1}\xi_{1} & 0\\
\end{array}\right ), \label{(3.4)}\\[-0.1cm]
\end{equation}
\begin{equation}
\overset{q}{\rho}_{1} = P_{1}^{\alpha}(x,k)=\left (
\begin{array}{ccc}
0 & k & -\alpha_{3} x^{2}\\[0.2cm]
-k & 0 & -\displaystyle\frac{\alpha_{2}\alpha_{3}}{\alpha_{1}} x^{1}\\[0.3cm]
\alpha_{3} x^{2} & \displaystyle\frac{\alpha_{2}\alpha_{3}}{\alpha_{1}} x^{1} & 0\\
\end{array}\right ),\label{(3.5)}
\end{equation}
\begin{equation}
\overset{q}{\rho}_{2} =-\left (
\begin{array}{ccc}
0 & (x^{3})^{q} & 0\\[0.3cm]
- (x^{3})^{q} & 0 & -\displaystyle\frac{\alpha_{3}}{\alpha_{1}} (x^{1})^{q}\\[0.3cm]
0 & \displaystyle\frac{\alpha_{3}}{\alpha_{1}} (x^{1})^{q} &  0\\
\end{array}\right )~~~\label{(3.6)}
\end{equation}
\begin{equation}
\begin{array}{lcl} H (x,\xi,k)& =&
\displaystyle\frac{1}{\Gamma(1+q)}[ - \alpha_{2}
(x^{1})^{q}(\xi_{1})^{q} + \alpha_{1} (x^{2})^{q}(\xi_{2})^{q} +
\displaystyle\frac{k \alpha_{1}}{\alpha_{3}}(\xi_{3})^{q}].
\label{(3.7)}
\end{array}
\end{equation}
\begin{Prop}
Let be $~k\neq 0 $ and $ q\in (0,1].$ The fractional system  $
(3.3) $ on the fractional Leibniz algebroid $ ({\bf R}^{3}\times
{\bf R}^{3}, \widetilde{P}, \overset{q}{\rho}_{1},
\overset{q}{\rho}_{2} ) $ associated to $ H $, where $
\widetilde{P},\overset{q}{\rho}_{1}, \overset{q}{\rho}_{2},
 H $ are given by the relations $ (3.4)-(3.7),$ is described by:
\begin{equation}
\left\{ \begin{array}{lcl} D_{t}^{q} \xi_{1} & = & k
\alpha_{1}((x^{2})^{q}- x^{2} \xi_{2}+
(\xi_{2})^{q})\\[0.2cm]
 D_{t}^{q} \xi_{2} & = & k \alpha_{2}((x^{1})^{q} -
x^{1} \xi_{1} + (\xi_{1})^{q})\\[0.2cm]
 D_{t}^{q} \xi_{3} &=&
\alpha_{2}\alpha_{3} (-(x^{1})^{q}x^{2}\xi_{2}+ x^{1}(x^{2})^{q}
\xi_{1} +
x^{1}(\xi_{2})^{q} - x^{2}(\xi_{1})^{q} ).\\[0.2cm]
D_{t}^{q} x^{1} & = & \alpha_{1}
(x^{2})^{q}(x^{3})^{q} \\[0.2cm]
 D_{t}^{q} x^{2} & = &
\alpha_{2} (x^{1})^{q}(x^{3})^{q} -k (x^{1})^{q}\\[0.2cm]
D_{t}^{q} x^{3} & = & \alpha_{3} (x^{1})^{q}(x^{2})^{q}
\label{(3.8)}
\end{array}\right.\\[-0.2cm]
\end{equation}
\end{Prop}
{\it Proof.} We have $~D_{\xi_{1}}^{q} H=- \alpha_{2}
(x^{1})^{q},~ D_{\xi_{2}}^{q} H= \alpha_{1} (x^{2})^{q},~
D_{\xi_{3}}^{q} H=\displaystyle\frac{k\alpha_{1}}{\alpha_{3}},~
D_{x^{1}}^{q} H = -\alpha_{2} (\xi_{1})^{q},~ D_{x^{2}}^{q} H=
\alpha_{1}(\xi_{2})^{q},~ D_{x^{3}}^{q} H= 0.$ Then $ (3.3) $ are
written as follows:
\[
\left (\begin{array}{c}
D_{t}^{q}{\xi}_{1}\\[0.2cm]
D_{t}^{q}{\xi}_{2}\\[0.2cm]
D_{t}^{q}{\xi}_{3}\\
\end{array}\right )= \widetilde{P} \left (\begin{array}{c}
- \alpha_{2}
(x^{1})^{q}\\[0.2cm]
\alpha_{1} (x^{2})^{q}\\[0.2cm]
 \displaystyle\frac{k\alpha_{1}}{\alpha_{3}}\\
\end{array}\right ) + \overset{q}{\rho}_{1} \left (\begin{array}{c}
-\alpha_{2} (\xi_{1})^{q}\\[0.2cm]
\alpha_{1} (\xi_{2})^{q} \\[0.2cm]
0 \\
\end{array}\right ),~~\left( \begin{array}{c}
D_{t}^{q} x^{1}\\[0.2cm]
D_{t}^{q} x^{2}\\[0.2cm]
D_{t}^{q} x^{3}\\
 \end{array}\right )= - \overset{q}{\rho}_{2} \left (\begin{array}{c}
- \alpha_{2}
(x^{1})^{q}\\[0.2cm]
\alpha_{1} (x^{2})^{q}\\[0.2cm]
 \displaystyle\frac{k\alpha_{1}}{\alpha_{3}} \\
\end{array}\right ).
\]

Replacing now  the matrices $ \widetilde{P}, \overset{q}{\rho}_{1}
$ and $ \overset{q}{\rho}_{2} $ one obtains the equations $ (3.8).
\hfill\Box $

The dynamical system $ (3.8)$ is called the  {\it $~q-$fractional
Euler top system with one control on the fractional Leibniz
algebroid} $ \pi : E = {\bf R}^{3}\times {\bf R}^{3}\to {\bf
R}^{3}$ for $ k\neq 0.$
\begin{Ex}
{\rm Taking $~k\neq 0, \alpha= (1, 1, 1) $ and $ q\in (0,1)$ in
the relations $ (3.4)-(3.8),$ the fractional system $ (3.3) $ on
the fractional Leibniz algebroid $ ({\bf R}^{3}\times {\bf R}^{3},
\widetilde{P}, \overset{q}{\rho}_{1}, \overset{q}{\rho}_{2} ) $
associated to function $ H,$
 is given by:
\begin{equation}
\left\{ \begin{array}{lcl} D_{t}^{q}\xi_{1} & = & k [(x^{2})^{q}-
x^{2} \xi_{2}+
(\xi_{2})^{q}]\\[0.2cm]
 D_{t}^{q}\xi_{2} & = & k [(x^{1})^{q} -
x^{1}\xi_{1}+ (\xi_{1})^{q} ]\\[0.2cm]
 D_{t}^{q}\xi_{3} &=&
-(x^{1})^{q}x^{2}\xi_{2}+ x^{1}(x^{2})^{q} \xi_{1} +
x^{1}(\xi_{2})^{q} - x^{2}(\xi_{1})^{q} \\[0.2cm]
D_{t}^{q} x^{1} & = &
(x^{2})^{q}(x^{3})^{q} \\[0.2cm]
 D_{t}^{q} x^{2} & = &
(x^{1})^{q}(x^{3})^{q} -k (x^{1})^{q}\\[0.2cm]
D_{t}^{q} x^{3} & = & (x^{1})^{q}(x^{2})^{q}.\label{(3.9)}
\end{array}\right.
\end{equation}

The dynamical system $ (3.9)$ is called the } {\it $~q-$fractional
Lagrange system with one control on the fractional Leibniz
algebroid $ \pi : E = {\bf R}^{3}\times {\bf R}^{3}\to {\bf
R}^{3}.\hfill\Box$}
\end{Ex}

\section{Numerical integration of the fractional Euler top with one control $ (3.8)$}

Let $f: {\bf R}\to {\bf R}$ be an integrable function. For $ q\in
(0, 1]$ we consider a function of
 $C^{1}-$class $~g^{q}: {\bf R}\times {\bf R} \to {\bf R},~(t,s)\to g_{t}^{q}(s)$
for $ t,s \in {\bf R}$ and $t_{0} \leq s \leq t.$

A {\it Riemann integral} of $f$ with respect to $ g_{t}^{q}$ is
defined by $~_{t_{0}} I_{t,g_{t}^{q}}f(t) =\int_{t_{0}}^{t} f(s)
g_{t}^{q}(s)ds.$

If the function $~g_{t}^{q}~$ is defined by\\[-0.2cm]
\begin{equation}
g_{t}^{q}(s) = \displaystyle\frac{1}{\Gamma(q)} (t-s)^{q-1}
e^{-\rho (t-s)},\label{(4.1)}
\end{equation}
where $\Gamma(q)$ is the Euler Gamma function and $\rho > 0,$ then
the Riemann integral $~_{t_{0}} I_{t,g_{t}^{q}}f(t)~$ becomes the
{\it fractional integral}, denoted with $_{t_{0}} I_{t}^{q}f(t)$,
  where:\\[-0.4cm]
\begin{equation}
_{t_{0}} I_{t}^{q}f(t) =
\displaystyle\frac{1}{\Gamma(q)}\int_{t_{0}}^{t}
(t-s)^{q-1}e^{-\rho (t-s)} f(s)ds.\label{(4.2)}
\end{equation}
\begin{Rem}
{\rm  For $\rho =0,$  the relation $(4.2)$ is the fractional
Riemann-Liouville integral (\cite{kilb}, Section 2.1). Also, the
fractional integral defined by $(4.2)$ is a special case of the
generalized fractional El-Nabulsi integral \cite{nabu}.}\hfill$\Box$
\end{Rem}
 If  $ \varphi:{\bf R}\times {\bf
R}^{n} \to {\bf R}^{n}$ is a deterministic function, we will
denote by:\\[-0.2cm]
\[
_{t_{0}} I_{t, g_{t}^{q}} \varphi(t,x(t)) =\int_{t_{0}}^{t}
\varphi(s,x(s)) g_{t}^{q}(s) ds,\\[-0.2cm]
\]
 where $~x(t)=(x^{1}(t),..., x^{n}(t)) $ and the integral is
Riemann integral  with respect to $ g_{t}^{q}.$

We call {\it Volterra integral equation with respect to $
g_{t}^{q},$} the functional Volterra type equation given by:\\[-0.4cm]
\begin{equation}
x(t) = x(t_{0}) +~ _{t_{0}}I_{t, g_{t}^{q}}
\varphi(t,x(t)).\label{(4.3)}
\end{equation}
The equation $(4.3)$ can be written formally in the following way:\\[-0.4cm]
\begin{equation}
dx(s) = \varphi(s,x(s))g_{t}^{q}(s)ds.\label{(4.4)}
\end{equation}
We apply the above considerations for the fractional Euler top system $(3.8)$\\[-0.2cm]
\begin{equation}
\left\{\begin{array}{lcl}
 D_{t}^{q} x^{i}(t) & = &
F_{i}(x^{1}(t), x^{2}(t), x^{3}(t)),~~~~~~~~~~~~~~~~~~~~~~~~~~~~ i=\overline{1,3}\\[0.2cm]
 D_{t}^{q} \xi_{i}(t) & = &
F_{i+3}(x^{1}(t), x^{2}(t), x^{3}(t), \xi_{1}(t), \xi_{2}(t), \xi_{3}(t)),~~~ i=\overline{1,3}\\[0.2cm]
x(t_{0}) &=& (x^{1}(t_{0}), x^{2}(t_{0}), x^{3}(t_{0})),~~~ \xi(t_{0}) = (\xi_{1}(t_{0}), \xi_{2}(t_{0}), \xi_{3}(t_{0}))
\end{array}\right.\label{(4.5)}\\
\end{equation}
where $~ t\in (t_{0}, \tau),~q \in (0,1)~$ and\\[-0.2cm]
\begin{equation}
\left\{\begin{array}{lcl}
F_{1}(x) & = & \alpha_{1}
(x^{2})^{q}(x^{3})^{q} \\[0.2cm]
F_{2}(x) & = & \alpha_{2} (x^{1})^{q}(x^{3})^{q} -k (x^{1})^{q}\\[0.2cm]
F_{3}(x) & = & \alpha_{3} (x^{1})^{q}(x^{2})^{q}\\[0.2cm]
F_{4}(x,\xi) & = & k
\alpha_{1}((x^{2})^{q}- x^{2} \xi_{2}+
(\xi_{2})^{q})\\[0.2cm]
F_{5}(x,\xi) & = & k \alpha_{2}((x^{1})^{q} -
x^{1} \xi_{1} + (\xi_{1})^{q})\\[0.2cm]
F_{6}(x,\xi) & = & \alpha_{2}\alpha_{3} (-(x^{1})^{q}x^{2}\xi_{2}+ x^{1}(x^{2})^{q}
\xi_{1} + x^{1}(\xi_{2})^{q} - x^{2}(\xi_{1})^{q} ).\label{(4.6)}
\end{array}\right.\\[-0.2cm]
\end{equation}
  Since the functions $ F_{i}(x)~$ and $~F_{i+3}(x,\xi)~$  for $ i=\overline{1,3} $ are continuous, the initial value problem $(4.5)$ is equivalent to system of Volterra integral equations \cite{kilb, podl}, which is given as follows:\\[-0.4cm]
\begin{equation}
\left\{\begin{array}{lcl}
x^{i}(t) & = & x^{i}(t_{0})  + ~_{t_{0}}I_{t, g_{t}^{q}}
F_{i}(x^{1}(t), x^{2}(t), x^{3}(t)),
~~~~~i=\overline{1,3},\\[0.2cm]
\xi_{i}(t) & =& \xi_{i}(t_{0})  + ~_{t_{0}}I_{t, g_{t}^{q}}
F_{i+3}(x^{1}(t), x^{2}(t), x^{3}(t), \xi_{1}(t), \xi_{2}(t), \xi_{3}(t)),
~~~~~i=\overline{1,3}.\label{(4.7)}
\end{array}\right.\\[-0.2cm]
\end{equation}
The equations $(4.7)$ can be written in the following form:\\[-0.2cm]
\begin{equation}
\left\{\begin{array}{lcl}
 dx^{1}(s) & = & F_{1}(x(s)) g_{t}^{q}(s)ds\\[0.2cm]
dx^{2}(s) & = & F_{2}(x(s)) g_{t}^{q}(s)ds,\\[0.2cm]
dx^{3}(s) & = & F_{3}(x(s))g_{t}^{q}(s)ds \\[0.2cm]
d\xi_{1}(s) & = & F_{4}(x(s), \xi(s)) g_{t}^{q}(s)ds\\[0.2cm]
d\xi_{2}(s) & = & F_{5}(x(s), \xi(s)) g_{t}^{q}(s)ds\\[0.2cm]
d\xi_{3}(s) & = & F_{6}(x(s), \xi(s)) g_{t}^{q}(s)ds,\label{(4.8)}
\end{array}\right.\\[-0.2cm]
\end{equation}
where $g_{t}^{q}~$ and
$~F_{i}(x(s)), ~ F_{i+3}(x(s), \xi (s)),~i=\overline{1,3}~ $  are given in $(4.1)$ and  $(4.6).$

The system $(4.8)$ is called the {\it Volterra integral equations
associated to Euler top system} $(4.6)$.

The problem for solving the system $(4.8)$ is reduced to one of
solving a sequence of systems of fractional equations in
increasing dimension on successive intervals $[j, (j+1)]$.

For the numerical integration of the system $(4.8)$ one can use
the first Euler scheme \cite{die2, kilb}, which is expressed as follows:\\[-0.4cm]
\begin{equation}
\left\{\begin{array}{lcl}
x^{i}(j+1) & = & x^{i}(j)+ h F_{i}(x^{1}(j), x^{2}(j),
x^{3}(j))g_{t}^{q}(j),~~~~~~~~~~~~~~~~~~~~~~~~~~~~ i=\overline{1,3},\\[0.2cm]
\xi_{i}(j+1) & =& \xi_{i}(j)  + h F_{i+3}(x^{1}(j), x^{2}(j), x^{3}(j), \xi_{1}(j), \xi_{2}(j), \xi_{3}(j))
g_{t}^{q}(j),~~~ i=\overline{1,3},\label{(4.9)}
\end{array}\right.\\[-0.2cm]
\end{equation}
where $ j=0,1,2,...,N,  h=\displaystyle\frac{T}{N}, T>0, N>0.~$

More precisely, if $t_{0}=0, \rho>0 $ and $ \varepsilon > 0, $ the numerical integration of the system $(4.9)$ is given by:\\[-0.4cm]
\begin{equation}
\left \{ \begin{array}{ll}
x^{1}(j+1) &= x^{1}(j)+
h~\displaystyle\frac{1}{\Gamma(q)}(t-j)^{q-1}\varepsilon^{-\rho
(t-j)}(\alpha_{1}
(x^{2})^{q}(j)(x^{3})^{q}(j) )\\[0.3cm]
x^{2}(j+1) &= x^{2}(j)+
h~\displaystyle\frac{1}{\Gamma(q)}(t-j)^{q-1}\varepsilon^{-\rho
(t-j)}(\alpha_{2} (x^{1})^{q}(j)(x^{3})^{q}(j) -k (x^{1})^{q}(j) )\\[0.3cm]
 x^{3}(j+1) &= x^{3}(j) +
h~\displaystyle\frac{1}{\Gamma(q)}(t-j)^{q-1}\varepsilon^{-\rho
(t-j)}( \alpha_{3} (x^{1})^{q}(j)(x^{2})^{q}(j) )\\[0.3cm]
\xi_{1}(j+1) &= \xi_{1}(j)+
h~\displaystyle\frac{1}{\Gamma(q)}(t-j)^{q-1}\varepsilon^{-\rho
(t-j)}(k \alpha_{1}((x^{2})^{q}(j)- x^{2}(j) \xi_{2}(j)+
(\xi_{2})^{q}(j)))\\[0.3cm]
\xi_{2}(j+1) &= \xi_{2}(j)+
h~\displaystyle\frac{1}{\Gamma(q)}(t-j)^{q-1}\varepsilon^{-\rho
(t-j)}(k \alpha_{2}((x^{1})^{q}(j) -
x^{1}(j) \xi_{1}(j) + (\xi_{1})^{q}(j)))\\[0.3cm]
 \xi_{3}(j+1) &= \xi_{3}(j) +
h~\displaystyle\frac{1}{\Gamma(q)}(t-j)^{q-1}\varepsilon^{-\rho
(t-j)}( \alpha_{2}\alpha_{3} (-(x^{1})^{q}(j) x^{2}(j)\xi_{2}(j)+ \\[0.3cm]
& + x^{1}(j)(x^{2})^{q}(j)
\xi_{1}(j) + x^{1}(j)(\xi_{2})^{q}(j) - x^{2}(j)(\xi_{1})^{q}(j)))\\[0.3cm]
x^{i}(0)&= \xi_{i}(0)= x_{e}^{i}+\varepsilon,~~~i=\overline{1,3}.\\[0.3cm]
\end{array}\right. \label{(4.10)}
\end{equation}
Using \cite{bora, die2}, we have that the numerical algorithm given by
$(4.10)$ is convergent.
\begin{Ex}
{\rm  Let us we present the numerical simulation of solutions of $q-$fractional
Lagrange system with one control on a fractional Leibniz algebroid $(3.9).$  Then, in $(4.10)$ we take: $~ \alpha_{1} = \alpha_{2} = \alpha_{3} = 1, k = 1 ~$ and $~q= 0.5. $
For to integrate this $~q-$fractional model, we consider $~h = 0.01,
\varepsilon= 0.01,  N = 500, t = 502 $ and the initial conditions
$~x^{1}(0)={\xi}_{1}(0)=\varepsilon,
  x^{2}(0)={\xi}_{2}(0)=\varepsilon, x^{3}(0)={\xi}_{3}(0)=\varepsilon + 0.5.$

Using the software Maple, the orbits $(n, x^{i}(n)),
i=\overline{1,3}~$ and $~ (n, {\xi}_{i}(n)), i=\overline{1,3}$ of
the fractional model $(3.9)$ are represented in Figures
1-3  and Figures 4-6.

In the coordinate systems $ Ox^{1}x^{2}x^{3} $ and $
O{\xi}_{1}{\xi}_{2}{\xi}_{3} ,$ the orbits $ (x^{1}(n), x^{2}(n)),
x^{3}(n)) $ and $ ( {\xi}_{1}(n)), {\xi}_{2}(n), {\xi}_{3}(n)) $
for the solutions of equations $(3.9)$ are represented in the
Figures 7-8.
\begin{center}
\begin{tabular}{ccc}
\includegraphics[width=5cm]{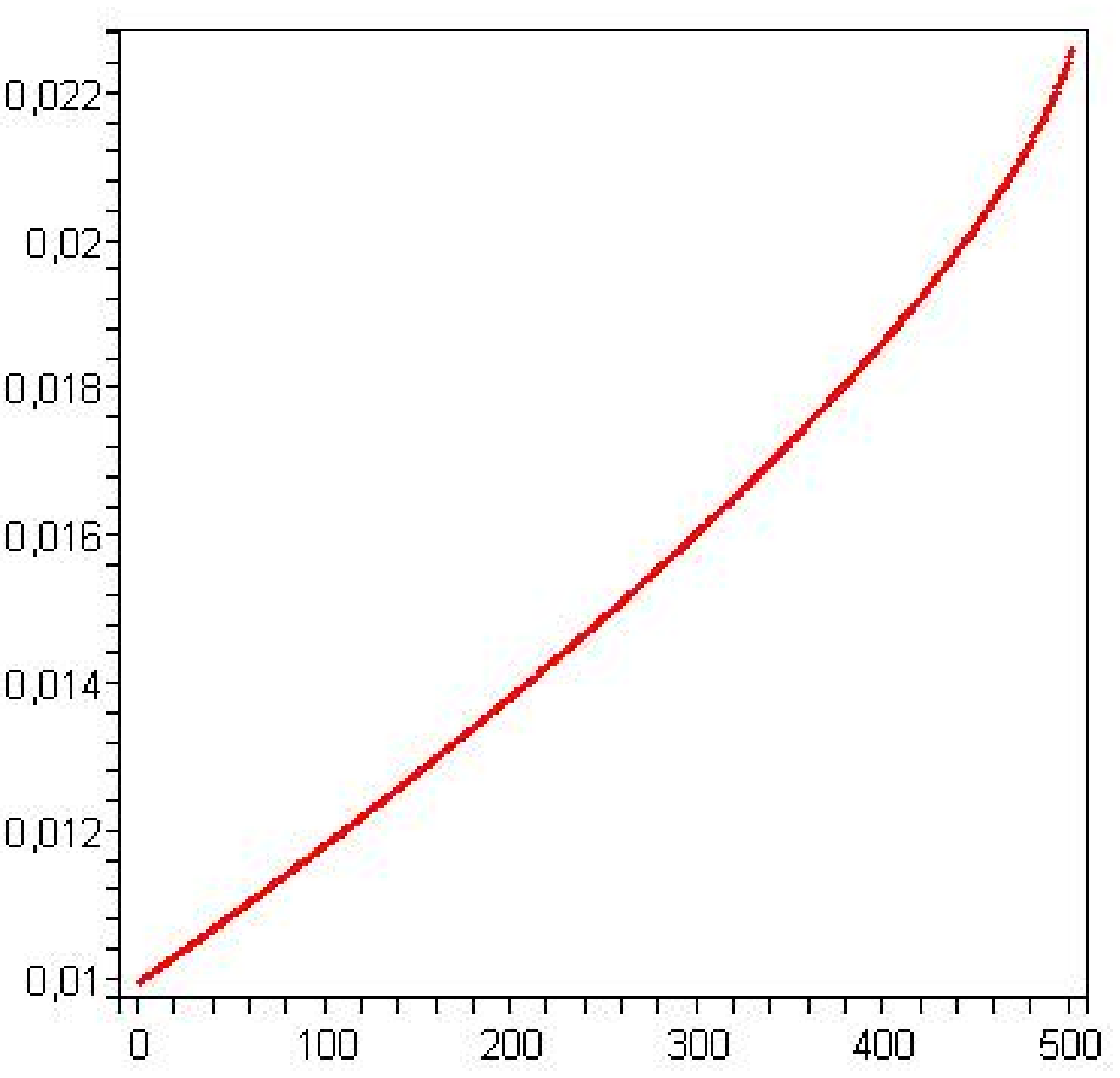}&
\includegraphics[width=5cm]{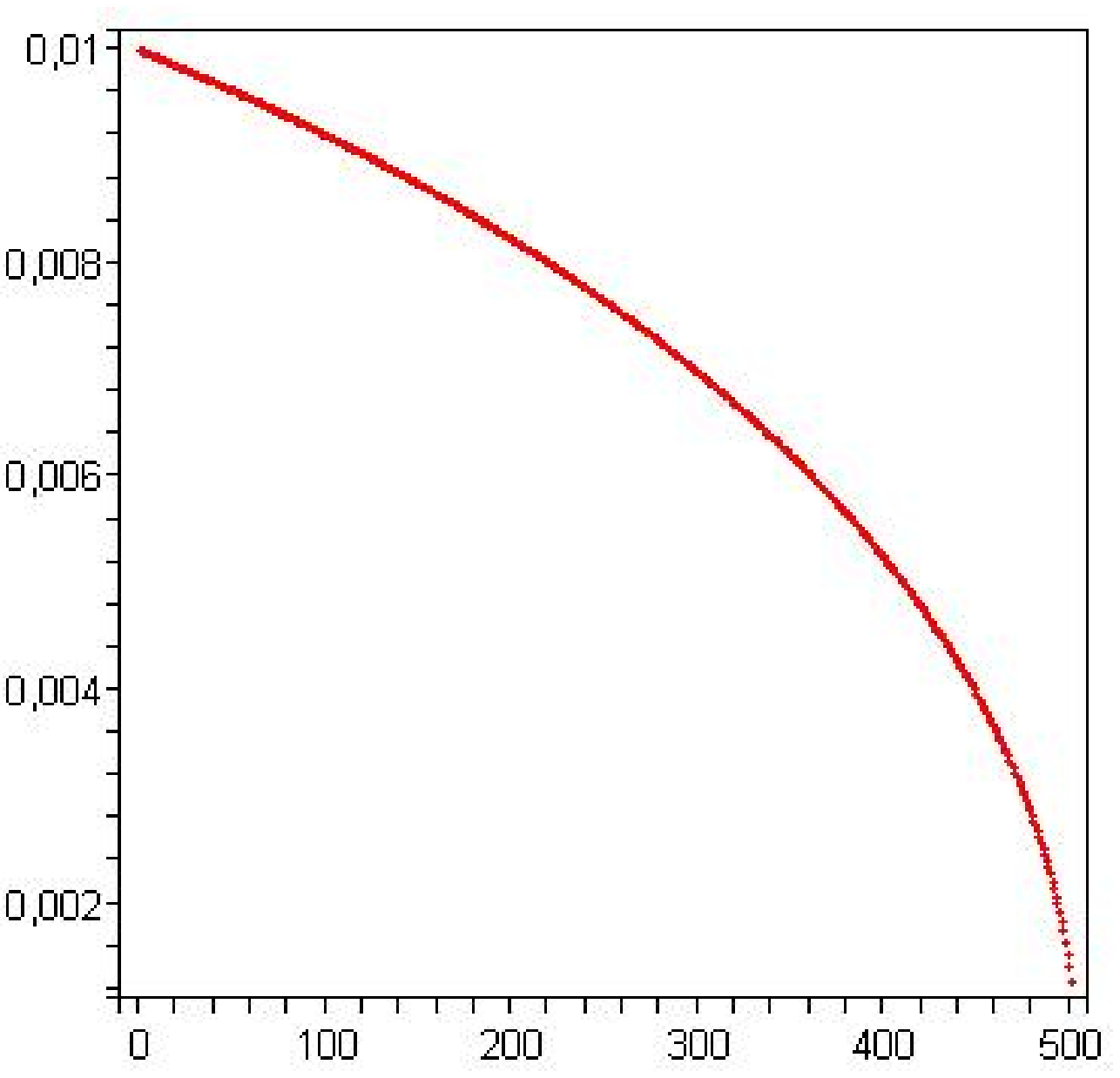}&
\includegraphics[width=5cm]{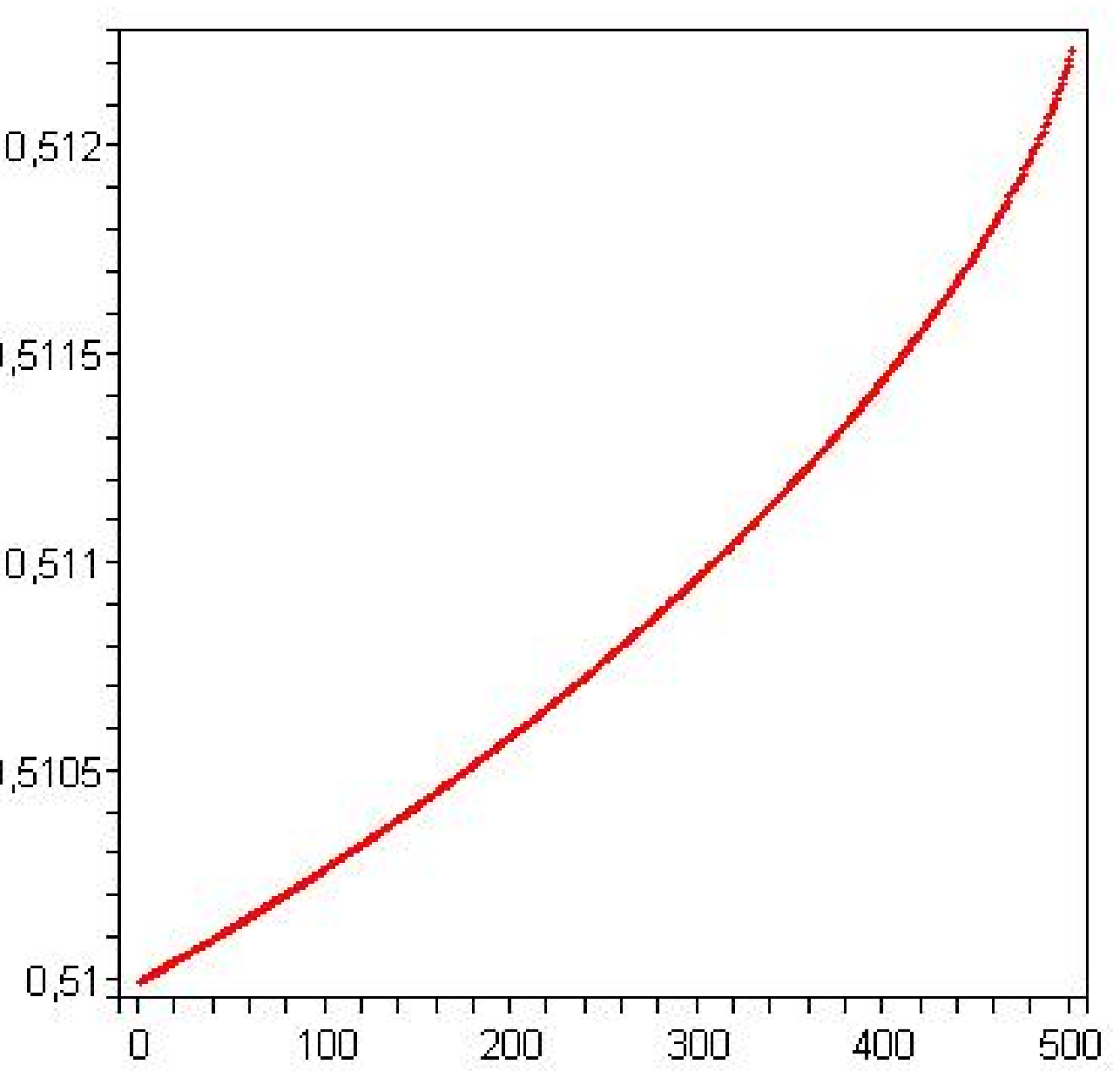}
\end{tabular}
\end{center}
\begin{center}
\begin{tabular}{ccc}
{\bf Figure 1.} $~(n, x^{1}(n))~~~~~~~~~~~~~~$ & {\bf Figure 2.} $~(n,
x^{2}(n))~~~~~~~~~~~~~~$ & {\bf Figure 3.} $~(n,
x^{3}(n))$\\
\end{tabular}
\end{center}
\begin{center}
\begin{tabular}{ccc}
\includegraphics[width=5cm]{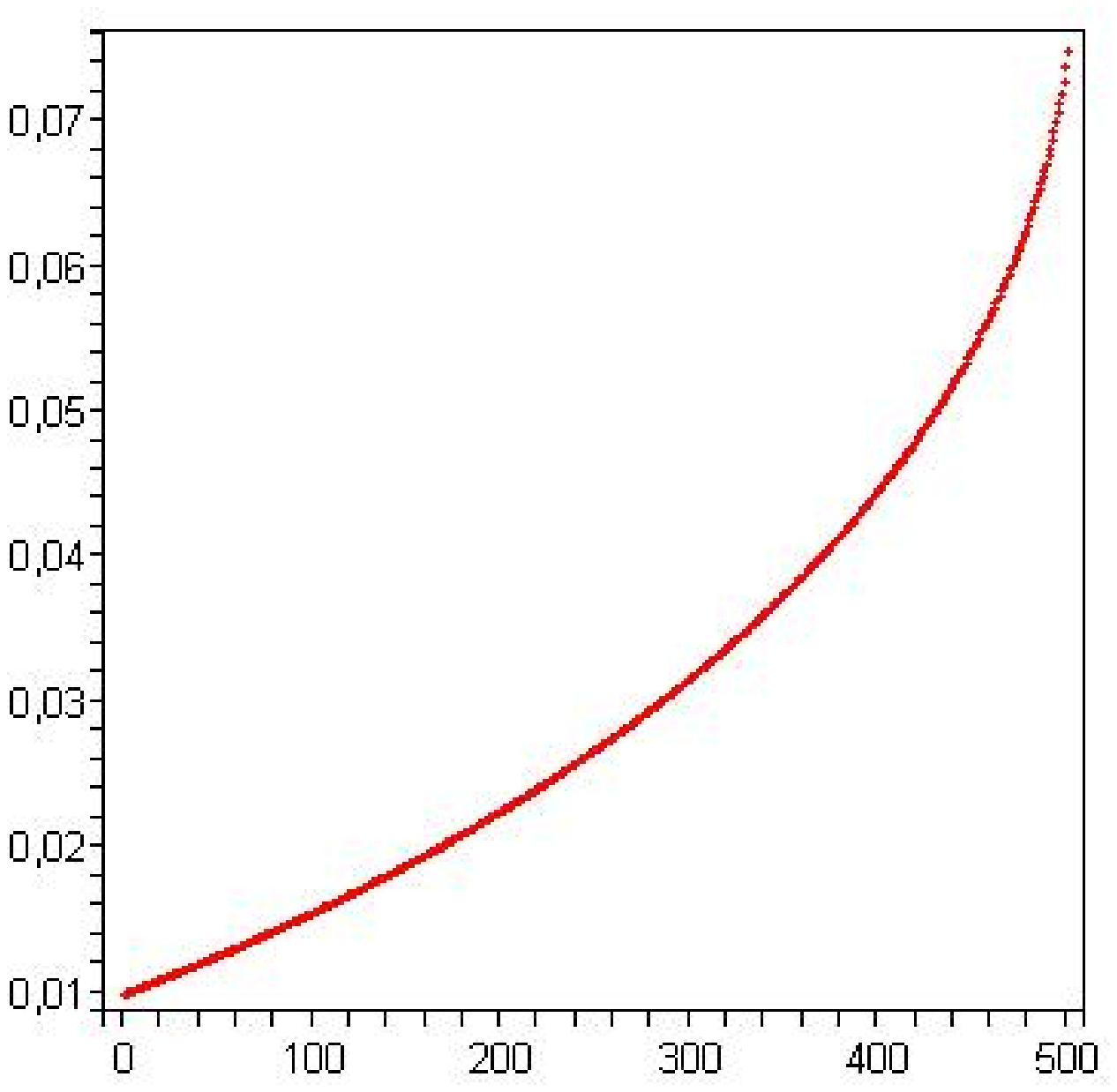}&
\includegraphics[width=5cm]{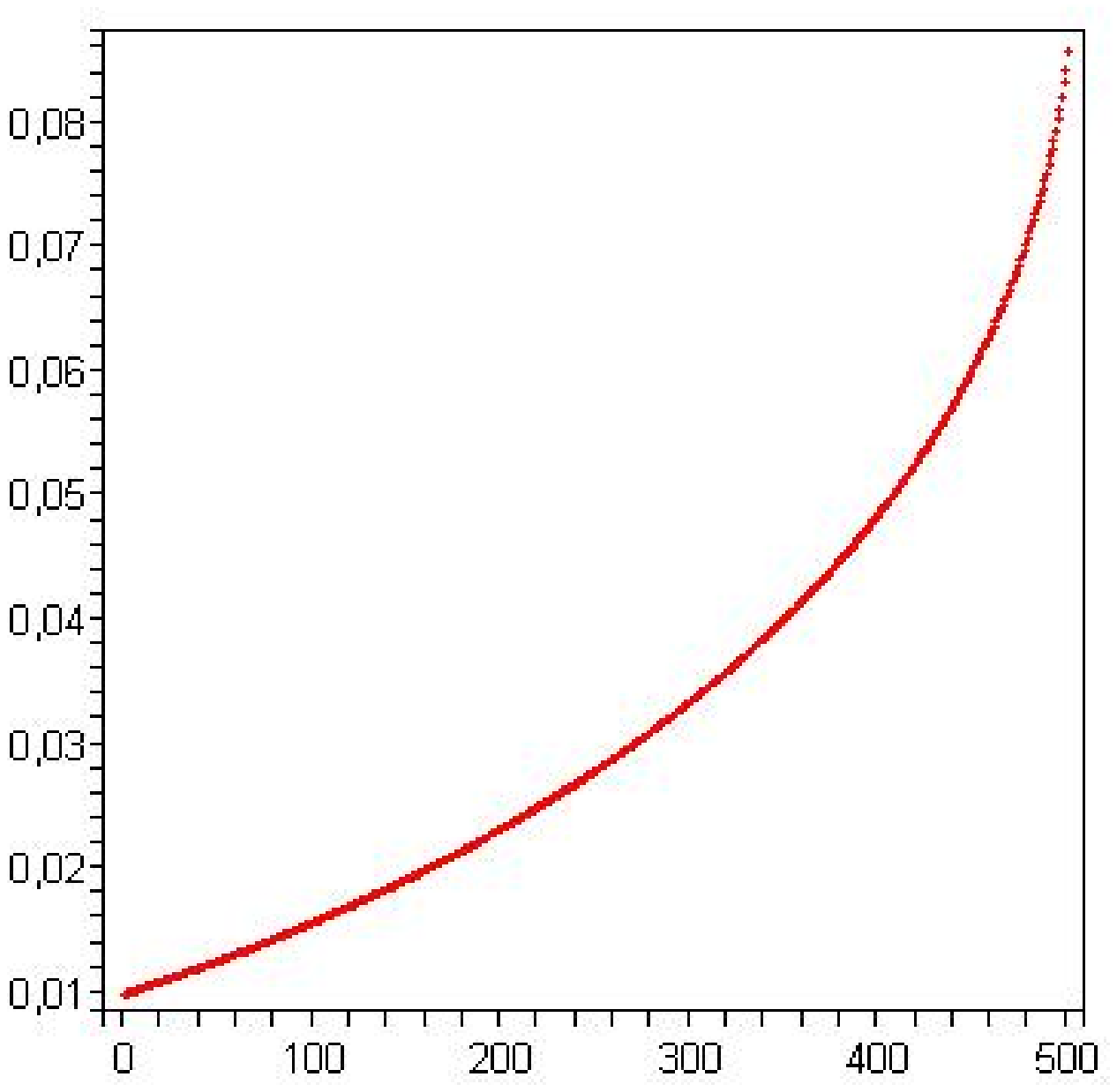}&
\includegraphics[width=5cm]{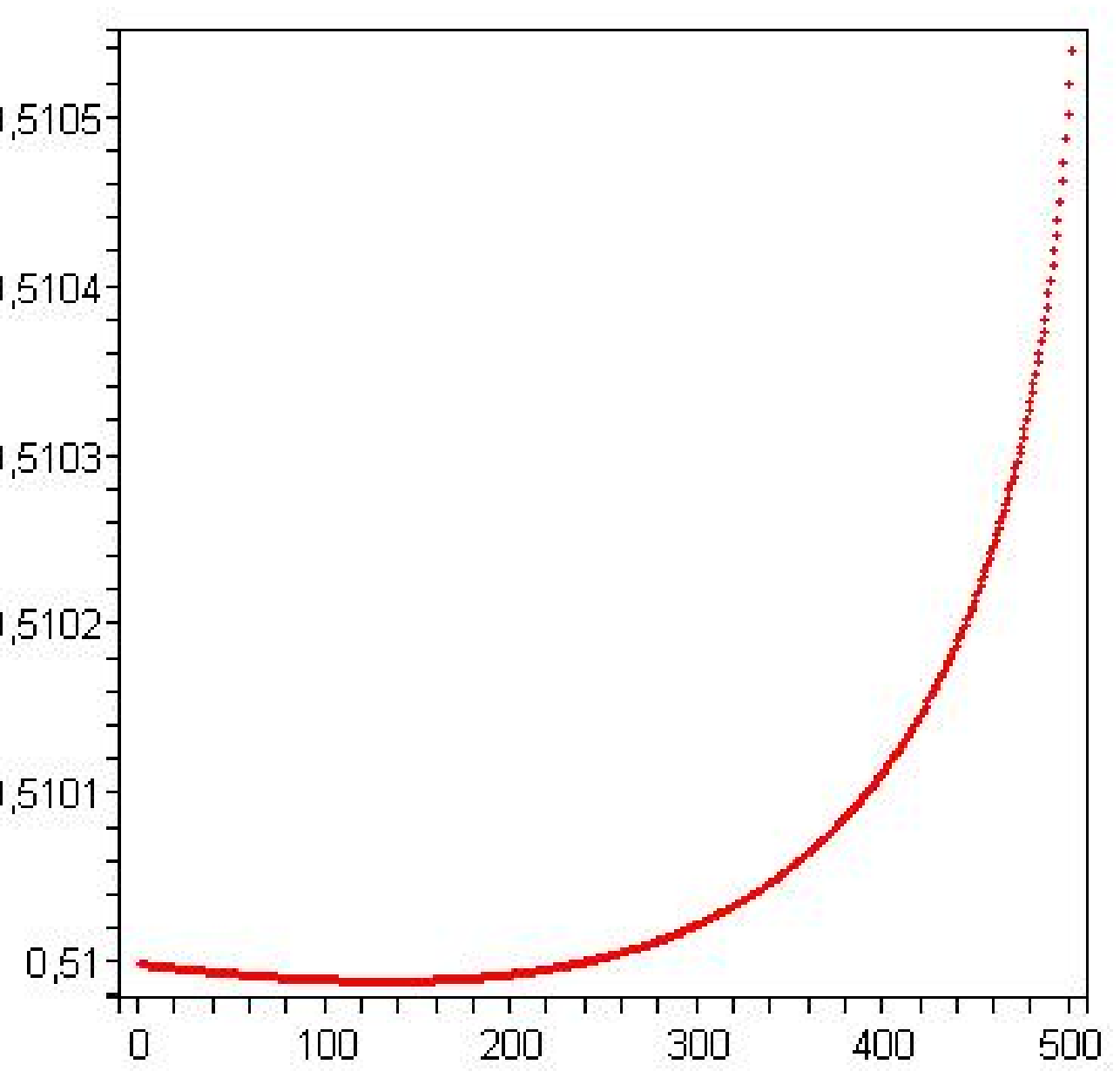}
\end{tabular}
\end{center}
\begin{center}
\begin{tabular}{ccc}
{\bf Figure 4.} $~(n, {\xi}_{1}(n))~~~~~~~~~~~~~~$ & {\bf Figure 5.}
$~(n, {\xi}_{2}(n))~~~~~~~~~~~~~~$ & {\bf Figure 6.} $~(n,
{\xi}_{3}(n))$\\
\end{tabular}
\end{center}
\begin{center}
\begin{tabular}{ccc}
\includegraphics[width=5cm]{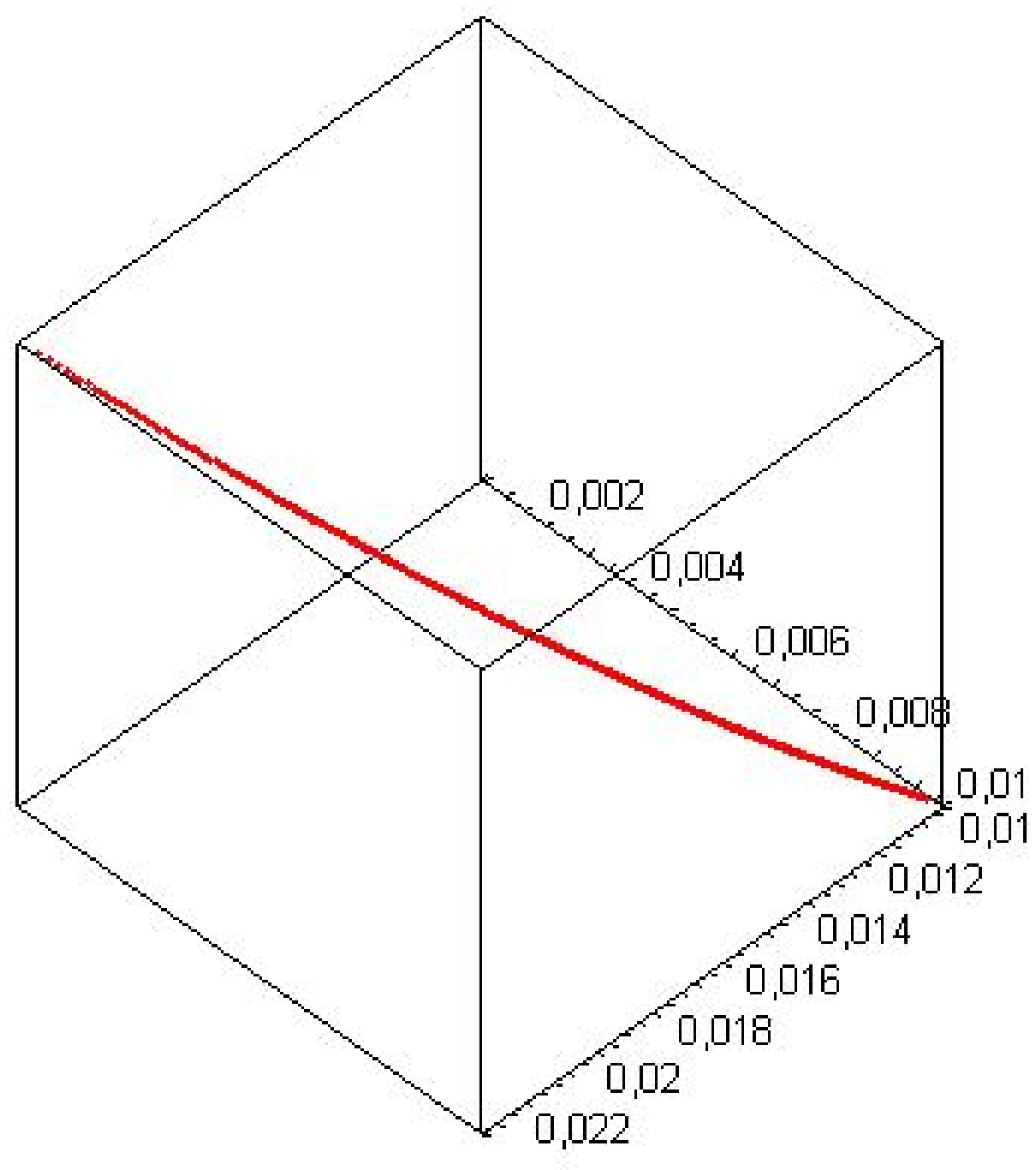} & &
\includegraphics[width=5cm]{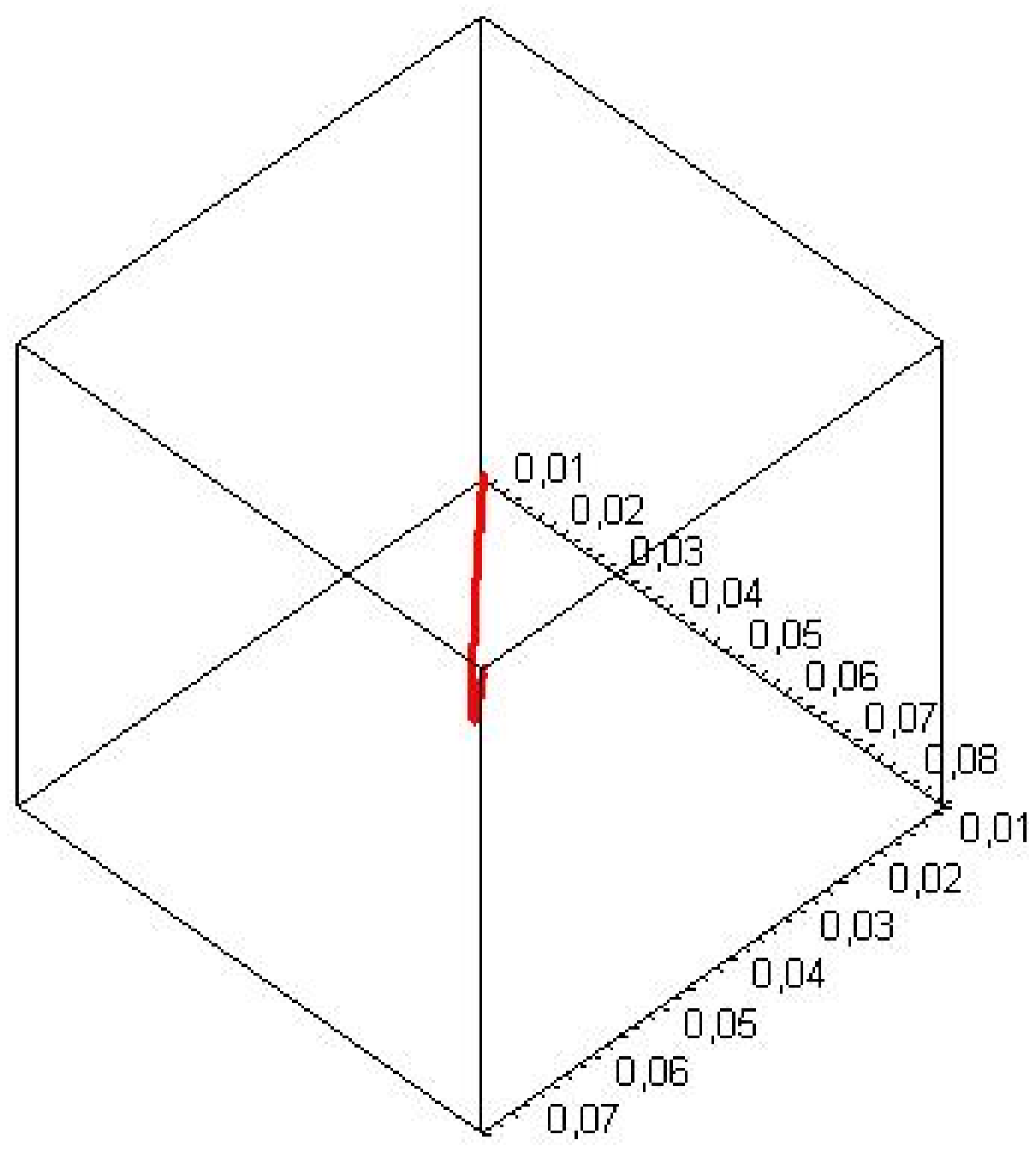}
\end{tabular}\\[-0.5cm]
\end{center}
\begin{center}
\begin{tabular}{ccc}
{\bf Figure 7.} $~(x^{1}(n), x^{2}(n), x^{3}(n)) $ for $ q=0.5$  & &
{\bf Figure 8.} $~(
{\xi}_{1}(n), {\xi}_{2}(n), {\xi}_{3}(n) ) $ for $ q=0.5.~$\hfill$\Box$
\end{tabular}
\end{center}}
\end{Ex}
 {\bf Conclusions.} This paper investigates  a new family of fractional differential systems, called the fractional Euler top system with one control $(2.4).$  The $ q-$fractional Euler top system with one control on a fractional Leibniz algebroid $(3.8) $ was constructed.  The  Adams-Bashforth-Moulton algorithm was applied for the numerical integration of the $ q-$fractional Euler top system $~(3.8).~$ Finally, the numerical simulation was  given for the $ q-$fractional Lagrange system with one control on a fractional Leibniz algebroid $(3.9).$ \hfill$\Box$\\[-0.8cm]

Author's adress\\[-0.2cm]

West University of Timi\c soara. Seminarul de Geometrie \c si Topologie.\\
Department of Mathematics. Timi\c soara, Romania.\\
E-mail: gheorghe.ivan@e-uvt.ro\\
\end{document}